\documentclass{amsart}
\usepackage{hyperref}

\newtheorem{thm}{Theorem}[section]
\newtheorem{lemma}[thm]{Lemma}
\newtheorem{coro}[thm]{Corollary}
\newtheorem{hypo}[thm]{Hypothesis {\bf H.}\hspace*{-0.6ex}}
\newtheorem{rem}[thm]{Remark}

\newcommand{\R}{{\mathbb R}}
\newcommand{\N}{{\mathbb N}}

\newcommand{\C}{{\mathbb C}}

\newcommand{\nn}{\nonumber}
\newcommand{\bea}{\begin{eqnarray}}
\newcommand{\eea}{\end{eqnarray}}
\newcommand{\ba}{\begin{array}}
\newcommand{\ea}{\end{array}}

\newcommand{\E}{\mathrm{e}}

\newcommand{\ti}{\tilde}
\newcommand{\bs}{\backslash}

\newcommand{\I}{\mathrm{i}}

\newcommand{\eps}{\varepsilon}
\newcommand{\sig}{\sigma}

\newcommand{\bth}{\begin{thm}}
\newcommand{\eth}{\end{thm}}
\newcommand{\bl}{\begin{lemma}}
\newcommand{\el}{\end{lemma}}
\newcommand{\bk}{\begin{coro}}
\newcommand{\ek}{\end{coro}}
\newcommand{\bh}{\begin{hypo}}
\newcommand{\eh}{\end{hypo}}
\newcommand{\br}{\begin{rem}}
\newcommand{\er}{\end{rem}}
\newcommand{\bpf}{\begin{proof}}
\newcommand{\epf}{\end{proof}}

\numberwithin{equation}{section}

\begin{document}

\title[On the Number of Eigenvalues of Jacobi Operators]{On the Finiteness of the
Number of Eigenvalues of Jacobi Operators below the Essential Spectrum}

\author{Franz Luef}
\address{Institut f\"ur Mathematik\\
Strudlhofgasse 4\\ 1090 Wien\\ Austria\\ and International Erwin Schr\"odinger
Institute for Mathematical Physics\\ Boltzmanngasse 9\\ 1090 Wien\\ Austria}

\author{Gerald Teschl}
\address{Institut f\"ur Mathematik\\
Strudlhofgasse 4\\ 1090 Wien\\ Austria\\ and International Erwin Schr\"odinger
Institute for Mathematical Physics\\ Boltzmanngasse 9\\ 1090 Wien\\ Austria}
\email{Gerald.Teschl@univie.ac.at}
\urladdr{http://www.mat.univie.ac.at/\string~gerald/}

\thanks{J. Difference Equ. Appl. {\bf 10}, no. 3, 299-307 (2004)}

%  Math Subject Classifications 

\maketitle

\begin{abstract}
We present a new oscillation criterion to determine whether 
the number of eigenvalues below the essential spectrum of
a given Jacobi operator is finite or not.  As an application we show
that Kenser's criterion for Jacobi operators follows as a special case.
\end{abstract}

\section{Introduction}

The goal of the present paper is to determine whether
the number of eigenvalues below the essential spectrum of the
Jacobi operator on $\ell^2(\N)$ associated with
\begin{equation}
(\tau f)(n) = a(n) f(n+1) + a(n-1) f(n-1) -b(n) f(n),
\end{equation}
where
\begin{equation}
a(n) \in \R \bs\{ 0\}, \quad b(n) \in \R, \quad n \in \N,
\end{equation}
is finite or not.

We will assume $a(n)<0$ (which is no restriction by \cite{tjac}, Lemma~1.6).
One of the main cases of interest is $a(n)= -1$ and one usually starts
with the operator $H_0$ associated with $b_0(n)= 2$. The spectrum is given
by $\sig(H_0)=[0,4]$. In particular, there are no eigenvalues below the essential
spectrum. Perturbing $b_0$ we can add any finite number of eigenvalues even if our
perturbation is of compact support. However, the question is, can we at least determine
whether the number of eigenvalues is finite or not, by looking at the asymptotics of
the perturbation? Moreover, what is the precise asymptotics separating the two
cases?

The natural tool for investigating such questions is oscillation theory
since finiteness of the number of eigenvalues is equivalent to the operator
being nonoscillatory. This fact first appeared in \cite{gl}. The
precise relation between the number of eigenvalues and the number of {\em
nodes} was established only recently by one of us in \cite{tosc}.

In the case of Sturm-Liouville operators there is a famous theorem by Kneser
\cite{kn} which gives a simple and beautiful answer to this question, with
many subsequent extensions by others. The most recent one being by \cite{gu},
who give a unified result containing all previously known ones as special
cases.

Unfortunately, it is not possible to use the same proof in the discrete case.
To understand why, let us first review the proof of Kneser's theorem in the
Sturm-Liouville case. The key idea is that the equation
\begin{equation}
\tau_0= -\frac{d^2}{dx^2} + \frac{\mu}{x^2}
\end{equation}
is of Euler type. Hence it is explicitly solvable with a fundamental system given
by
\begin{equation}
x^{\frac{1}{2} \pm \sqrt{\mu+\frac{1}{4}}}.
\end{equation}
There are two cases to distinguish. If $\mu\ge - 1/4$ all solutions are 
nonoscillatory. If $\mu< - 1/4$ one has to take real/imaginary parts and
all solutions are oscillatory. Hence a straightforward application of Sturm's
comparison theorem between $\tau_0$ and
\begin{equation}
\tau= -\frac{d^2}{dx^2} + q(x)
\end{equation}
yields
\begin{equation}
\lim_{x \to\infty} \ba{c}\inf \\ \sup \ea
\big( x^2 q(x) \big) \ba{c}> \\ <\ea -\frac{1}{4}
\mbox{ implies }
\ba{c} \text{nonoscillation} \\ \text{oscillation} \ea
\mbox{ of $\tau$ near $\infty$}.
\end{equation}
Since Sturm's comparison theorem is also available for Jacobi operators (see,
e.g, \cite{tjac}, Lemma~4.4) it seems easy to generalize this result by
considering the discrete Euler equation
\begin{equation}
u(n+1) - 2 u(n) + u(n-1) - \frac{\mu}{n(n-1)} u(n-1)=0.
\end{equation}
However, unfortunately, this equation is not symmetric! The corresponding results
for this equation can be found as special cases in \cite{arg}, Section~6.11.
Thus a straightforward generalization is not possible.

One approach is to consider generalized ordinary differential expressions which
contain difference equations as a special case. This can be found in \cite{gl} and
\cite{min} (see Section~2.2). As an alternative, we will prove a new
oscillation criterion and show that Kneser's criterion follows as a special case.

For further (non-)oscillation criteria we refer to Hinton and Lewis \cite{hl}
and Hooker and Patula \cite{hp} (see also \cite{har}, \cite{hkp}, \cite{pat1},
and \cite{pat2}). For a different approach using Birman-Schwinger type arguments see \cite{ger2},
\cite{ger3}.

Our present paper was motivated by the work of Gesztesy and \"Unal \cite{gu} mentioned earlier. In fact, it can be viewed as a discrete generalization of their
results. However, again a straightforward generalization is not possible since
their proofs also rely on explicit solubility of the involved equations.

\section{Main results and applications}

Before we can write down our main result, we need to fix some notation.
Recall that $\tau$ is called oscillatory if one (and hence any) real-valued
solution of $\tau u =0$ has an infinite number of nodes, that is, points $n \in
\N$, such that either
\begin{equation}
u(n) =0 \quad \mbox{or} \quad a(n)u(n)u(n+1)>0.
\end{equation}
In the special case $a(n)<0$, $n \in \N$, a node of $u$ is precisely a sign flip
of $u$ as one would expect. In the general case, however, one has to take the sign
of $a(n)$ into account. 

Recall that if $u_0(n)>0$ solves
\begin{equation} \label{deftaun}
(\tau_0 u)(n) = a(n) u_0(n+1) + a(n-1) u_0(n-1) -b_0(n) u_0(n) =0,
\end{equation}
then
\begin{equation}
\hat{u}_0(n)= u_0(n) Q_0(n), \quad Q_0(n)= \sum_{j=0}^{n-1} \frac{-1}{a(j)u_0(j)u_0(j+1)},
\end{equation}
is a second, linearly independent positive solution. A positive solution is
called minimal if
\begin{equation} \label{condminimal}
\lim_{n\to\infty} Q_0(n) = \infty.
\end{equation}
Minimal solutions are unique up to a multiple. See \cite{tjac}, Section~2.3 for more
information. With this notation our main result reads as follows:

\bth \label{thmmain}
Suppose $a(n)\in\R\bs\{0\}$, $b(n), b_0(n)\in\R$ such that $a_0<|a(n)|<A_0$
for some real constants $0<a_0<A_0$ and let $u_0$ be a non-decreasing minimal
positive solution of $\tau_0 u_0 =0$ (as defined in (\ref{deftaun})).
Introduce
\begin{equation}
A(n) = \frac{2 a(n-1) a(n+1)}{a(n-1)+a(n+1)}, \quad n\in\N,
\end{equation}
and
\begin{equation}
Q_0(n)= \sum_{j=0}^{n-1} \frac{-1}{a(j)u_0(j)u_0(j+1)}, \quad n\in\N.
\end{equation}
Then $\tau$ is nonoscillatory if
\begin{equation}
\liminf_{n\to\infty} -A(n) u_0^4(n) Q_0^2(n) (b(n)-b_0(n)) > -\frac{1}{4}
\end{equation}
and oscillatory if
\begin{equation}
\limsup_{n\to\infty} -A(n) u_0^4(n) Q_0^2(n) (b(n)-b_0(n)) < -\frac{1}{4}.
\end{equation}
\eth

The proof will be given in Section~\ref{secproof} below.

As a first application, let us show how this result can be used to answer our
question posed in the introduction. We choose
\begin{equation}
a(n)=-1, \quad b_0(n)=2.
\end{equation}
Then we have
\begin{equation}
u_0(n)=1 \quad\mbox{and}\quad \hat{u}_0(n)=n
\end{equation}
and thus
\begin{equation}
\lim_{n \to\infty} \ba{c}\inf \\ \sup \ea
\big( n^2 (b(n)-2) \big) \ba{c}> \\ <\ea -\frac{1}{4}
\mbox{ implies }
\ba{c} \text{nonoscillation} \\ \text{oscillation} \ea
\mbox{ of $\tau$ near $\infty$},
\end{equation}
which is the claimed generalization of Kneser's result. Clearly, the next
question is what happens in the limiting case, where
$\lim_{n\to\infty} n^2(b(n)-2)= -4^{-1}$? This can be answered by our result
as well:

Recall the iterated logarithm
\begin{equation}
\ln_0(x)=x, \qquad \ln_k(x)=\ln_{k-1}(\ln(x)),
\end{equation}
where $\ln_k(x)$ is defined for $x>\E_k$, with $\E_1=0$, $\E_k=\E^{\E_{k-1}}$.

\bk
Let 
\begin{equation}
a(n)=-1, \qquad b_k(n)= 2 - \frac{1}{4} \sum_{j=0}^{k-1} \frac{1}{\prod_{\ell=0}^j
\ln_\ell(n)^2}.
\end{equation}
Then $\tau$ is nonoscillatory if 
\begin{equation}
\liminf_{n\to\infty} \big( \prod_{j=0}^k \ln_j(n) \big)^2 (b(n)-b_k(n))  >
-\frac{1}{4}
\end{equation}
and oscillatory if
\begin{equation}
\limsup_{n\to\infty} \big( \prod_{j=0}^k \ln_j(n) \big)^2 (b(n)-b_k(n))  <
-\frac{1}{4}.
\end{equation}
\ek

\bpf
To show how this follows from our result we consider
\begin{equation}
u_k(n)= \sqrt{\prod_{j=0}^{k-1} \ln_j(n)},
\end{equation}
which is a solution of $\ti{\tau}_k$ associated with
\begin{equation}
a(n)= -1, \qquad \ti{b}_k(n)= \frac{u_k(n+1) + u_k(n-1)}{u_k(n)}.
\end{equation}
To prove the claim it suffices to show
\begin{eqnarray} \nn
\ti{b}_k(n) &=& b_k(n) + O(n^{-3}),\\ \label{diffbtb}
Q_k(n) &=& \ln_k(n) + O(1)
\end{eqnarray}
since the differences will not contribute to the limits from above.

To establish (\ref{diffbtb}) we first recall the following formulas for the
first and second derivative of $\ln_k(x)$:
\begin{eqnarray} \nn
\ln_k'(x) &=& \prod_{j=0}^{k-1} \frac{1}{\ln_j(x)},\\ 
\ln_k''(x) &=&  - \ln_k'(x) \sum_{j=1}^k \ln_j'(x), \qquad x>\E_k.
\end{eqnarray}
Now we can show (\ref{diffbtb}). First of all we have
\begin{equation}
Q_k(n)= \int^n \frac{dx}{u_k(x)^2} + O(1) = \ln_k(n) +O(1).
\end{equation}
The second claim is a bit harder. We begin with
\begin{eqnarray} \nn
\frac{u_k(n\pm 1)}{u_k(n)} &=& \left(\prod_{j=0}^{k-1} \frac{\ln_j(n\pm 1)}{\ln_j(n)}
\right)^{1/2}= \prod_{j=0}^{k-1} \left( \sum_{\ell=0}^\infty
\frac{(-1)^\ell}{\ell!} \frac{\ln_j^{(\ell)}(n)}{\ln_j(n)} \right)^{1/2}\\ \nn
&=& \prod_{j=0}^{k-1} \left( 1 \pm \frac{\ln_j'(n)}{\ln_j(n)} + \frac{1}{2} 
\frac{\ln_j''(n)}{\ln_j(n)} + O(n^{-3}) \right)^{1/2}\\ \nn
&=& \prod_{j=0}^{k-1} \left( 1 \pm \frac{1}{2} \ln_{j+1}'(n) + \frac{1}{4} 
\left( \frac{\ln_j''(n)}{\ln_j(n)} - \frac{\ln_{j+1}'(n)^2}{2} \right) + O(n^{-3})
\right)\\ \nn
&=& 1 \pm \frac{1}{2} \sum_{j=0}^{k-1} \ln_{j+1}'(n) + \frac{1}{4} \sum_{j=0}^{k-1}
\left( \frac{\ln_j''(n)}{\ln_j(n)} - \frac{\ln_{j+1}'(n)^2}{2} \right) +\\ \nn
&& {} + \frac{1}{4} \sum_{j=0}^{k-1} \ln_{j+1}'(n) \sum_{\ell=0}^{j-1} \ln_{\ell+1}'(n) +
O(n^{-3}). \\
&=& 1 \pm \frac{1}{2} \sum_{j=0}^{k-1} \ln_{j+1}'(n) - \frac{1}{8} \sum_{j=0}^{k-1}
\ln_{j+1}'(n)^2 + O(n^{-3}).
\end{eqnarray}
Now combining both formulas we obtain the desired result
\begin{equation}
\ti{b}_k(n) = 2 - \frac{1}{4} \sum_{j=0}^{k-1} \ln_{j+1}'(n)^2 + O(n^{-3})\\
= b_k(n) + O(n^{-3}).
\end{equation}
\epf

Another interesting example is the case
\begin{equation}
b_0(n)= - a(n) - a(n-1).
\end{equation}
Again we can take $u_0(n)=1$ to obtain

\bk
Let $a_0 \le |a(n)|\le A_0$ and abbreviate
\begin{equation}
A(n) = \frac{2 a(n-1) a(n+1)}{a(n-1)+a(n+1)}, \qquad
Q_0(n) = \sum_{j=0}^{n-1} \frac{-1}{a(j)}.
\end{equation}
Then $\tau$ is nonoscillatory if 
\begin{equation}
\liminf_{n\to\infty} -A(n) Q_0(n)^2 (b(n) +a(n-1)+a(n)) > -\frac{1}{4}
\end{equation}
and oscillatory if
\begin{equation}
\limsup_{n\to\infty} -A(n) Q_0(n)^2 (b(n) +a(n-1)+a(n)) < -\frac{1}{4}.
\end{equation}
\ek

Of course one could take two arbitrary sequence $a(n)<0$ and $u_0(n)>0$ such
that $u_0$ is non-decreasing and (\ref{condminimal}) is satisfied,
compute $b_0(n) = -(a(n) u_0(n+1) + a(n-1) u_0(n-1))/u_0(n)$, and apply
Theorem~\ref{thmmain} to get a new (non)oscillation criterion.

\section{Proof of Theorem~\ref{thmmain}}
\label{secproof}

Suppose $a(n)\in\R\bs\{0\}$, $b(n), b_0(n)\in\R$ such that
\begin{equation}
a_0 \le |a(n)|\le A_0
\end{equation}
for some real constants $0<a_0<A_0$ and let $u_0$ be a non-decreasing minimal
positive solution of $\tau_0 u_0 =0$ as in the previous section. Note that the corresponding second positive solution $\hat{u}_0(n)$ is increasing.

First we collect some basic facts which will be needed later on.

\bl \label{lempopun}
Let $u_0$ be a minimal positive non-decreasing solution, then
we have
\begin{equation}
\lim_{n\to\infty} \frac{u_0(n+1)}{u_0(n)} = \lim_{n\to\infty} \frac{u_0(n)}{u_0(n-1)}
= 1
\end{equation}
and
\begin{equation}
u_0(n) \hat{u}_0(n) = u_0^2(n) Q_0(n) \ge \frac{n}{A_0}.
\end{equation}
\el

\bpf
Monotonicity of $u_0$ implies
\begin{equation}
\frac{1}{u_0(j+1)^2} \le \frac{1}{u_0(j)u_0(j+1)} \le \frac{1}{u_0(j)^2}.
\end{equation}
Summing the last expression from $0$ to $n-1$ and subtracting the right side
yields
\begin{equation}
0\le\sum_{j=0}^{n-1}\frac{1}{u_0(j)u_0(j+1)}\biggr(\frac{u_0(j+1)}{u_0(j)}-1\biggl)\le
\frac{1}{u_0(0)^2} - \frac{1}{u_0(n)^2} \le \frac{1}{u_0(0)^2}.
\end{equation}
Since $u(n)$ is minimal,
\begin{equation}
\sum_{j=0}^{n-1}\frac{1}{u_0(j)u_0(j+1)} \ge a_0 Q_0(n) \to \infty
\end{equation}
implies the first result.

For the second claim we use
\begin{equation}
Q_0(n) = \sum_{j=0}^{n-1}\frac{-1}{a(j) u_0(j)u_0(j+1)} \ge \frac{1}{A_0}
\sum_{j=0}^{n-1}\frac{1}{u_0(j+1)^2} \ge \frac{n}{A_0 u_0(n)^2}
\end{equation}
finishing the proof.
\epf

Our next goal is to find a suitable comparison equation. We do this by trying the ansatz
\begin{equation}
u_1(n) = u_0(n) Q_0(n)^\alpha.
\end{equation}
Then $u_1(n)$ satisfies
\begin{equation}
\tau_1 u(n)= a(n) u_1(n+1) + a(n-1) u_1(n-1) + b_1(n) u_1(n) =0
\end{equation}
with $b_1(n)$ given by
\begin{eqnarray} \nn
b_1(n) &=& -\frac{a(n) u_0(n+1)}{u_0(n)} \left( 1 - \frac{1}{a(n)u_0(n+1)^2 Q_0(n)}
\right)^\alpha\\
&& {} - \frac{a(n-1) u_0(n-1)}{u_0(n)} \left( 1 + \frac{1}{a(n-1)u_0(n-1)^2 Q_0(n)}
\right)^\alpha.
\end{eqnarray}

In order to get an oscillating comparison equation we need to admit $\alpha\in\C$.
However, this will also render $b_1(n)$ complex and hence it will be of no use
for us. To overcome this problem we look at the asymptotic behavior of
$b_1(n)$ for $n\to\infty$, which is given by
\begin{equation}
b_1(n) = b_0(n) + \mu U(n) + O(\frac{1}{u_0^6(n) Q_0^3(n)}), 
\quad \mu =\alpha(\alpha-1),
\end{equation}
where
\begin{equation}
U(n)= \frac{1}{2 u_0^4(n) Q_0^2(n)} \left(
\frac{-u_0(n)}{a(n+1) u_0(n+1)} + \frac{-u_0(n)}{a(n-1)u_0(n-1)} \right).
\end{equation}
If $\alpha\in\R$ we can choose $b_1(n)$ directly as comparison potential to
obtain that $\tau$ is nonoscillatory if
\begin{equation}
\liminf_{n\to\infty} \frac{b(n)-b_0(n)}{b_1(n)-b_0(n)} > \mu.
\end{equation}
Using the optimal value $\alpha=\frac{1}{2}$ plus the expansion
from above we end up with
\begin{equation}
\liminf_{n\to\infty} \frac{b(n)-b_0(n)}{U(n)} > -\frac{1}{4}.
\end{equation}

This settles the first part of our theorem. Now we come to the harder one.
As already noticed, in order to get an oscillating comparison equation we
need to choose complex values for $\alpha$. Our strategy is to choose
$\alpha= \frac{1}{2} + \I\eps$ such that at least $\mu= -\frac{1}{4}-\eps^2$
remains real and take $\ti{b}_1(n)= b_0(n) + \mu U(n)$ as comparison equation. Of course
we do not know the solutions of this equation, but our hope is that they
are asymptotically given by the real/imaginary parts of
\begin{equation}
u_1(n) = u_0(n) \sqrt{Q_0(n)} \left( \cos(\eps \ln Q_0(n)) + \I
 \sin(\eps \ln Q_0(n)) \right).
\end{equation}
Hence if we can show that there are solutions $\ti{u}_1$ of $\ti{\tau}_1 \ti{u}_1=0$
satisfying
\begin{equation} \label{asympuo}
\ti{u}_1(n) = u_1(n) ( 1 + o(1))
\end{equation}
we are done.

To show this, we begin with
\begin{equation}
\tau_1 \ti{u}_1(n) = \Delta(n) \ti{u}_1(n), \qquad \Delta(n) = b_1(n)-\ti{b}_1(n),
\end{equation}
and use the solution formula for the inhomogeneous equation (\cite{tjac},
eqn. (1.48)) to obtain the following equation
\begin{equation}
\ti{u}_1(n) = u_1(n) - \sum_{j=n+1}^{n_0} u_1(n) u_1(j) (Q_1(n)-Q_1(j)) \Delta(j)
\ti{u}_1(j),
\end{equation}
where $Q_1$ is defined as $Q_0$ but with $u_1$ in place of $u_0$. Formally,
letting $n_0\to\infty$, and setting
\begin{equation}
\ti{u}_1(n) = u_1(n) v(n)
\end{equation}
we obtain
\begin{equation}
v(n) = 1 - \sum_{j=n+1}^\infty u_1(j)^2 (Q_1(n)-Q_1(j)) \Delta(j) v(j).
\end{equation}
If we can show existence of a solution $v(n)=1 + o(1)$ of this last
equation, we are done. For this it suffices to verify the
assumptions of \cite{tjac}, Lemma~7.8. Hence we need to estimate
the kernel of the above sum equation. Using
\begin{eqnarray} \nn
|Q_1(n)-Q_1(j)| & \le& \sum_{k=n}^j \frac{1}{a(k) u_0(k) u_0(k+1) \sqrt{Q_0(k) Q_0(k+1)}}\\
&\le& \frac{1}{a_0} \sum_{k=n}^j \frac{1}{u_0(k)^2 Q_0(k)}
\end{eqnarray}
and
\begin{equation}
|\Delta(j)| \le \frac{const}{u_0(j)^6 Q_0(j)^3}
\end{equation}
we obtain by Lemma~\ref{lempopun}
\begin{eqnarray} \nn
|u_1(j)^2 (Q_1(n)-Q_1(j)) \Delta(j)| &\le& \frac{const}{u_0(j)^4 Q_0(j)^2}
\sum_{k=n}^j \frac{1}{u_0(k)^2 Q_0(k)}\\
&\le& const \frac{\ln(j)}{j^2}.
\end{eqnarray}
Thus we can apply \cite{tjac}, Lemma~7.8 to conclude existence of a
solution of type (\ref{asympuo}) which finishes the proof of Theorem~\ref{thmmain}.

\end{document}